\documentclass[11pt]{amsart}
\usepackage{amsopn}
\usepackage{amssymb, amscd, amsmath}
\usepackage[labelsep=period,indention=1.8cm,labelfont=bf,
textfont=it,justification=centering]{caption}
\usepackage{enumerate}

\newcommand{\nc}{\newcommand}

\nc{\vg}{\mathfrak{v} } \nc{\wg}{\mathfrak{w} } \nc{\zg}{\mathfrak{z} }
\nc{\ngo}{\mathfrak{n} } \nc{\kg}{\mathfrak{k} } \nc{\mg}{\mathfrak{m} }
\nc{\bg}{\mathfrak{b} } \nc{\ggo}{\mathfrak{g} } \nc{\ggob}{\overline{\mathfrak{g}}
} \nc{\sog}{\mathfrak{so} } \nc{\sug}{\mathfrak{su} } \nc{\spg}{\mathfrak{sp} }
\nc{\slg}{\mathfrak{sl} } \nc{\glg}{\mathfrak{gl} } \nc{\cg}{\mathfrak{c} }
\nc{\rg}{\mathfrak{r} } \nc{\hg}{\mathfrak{h} } \nc{\tg}{\mathfrak{t} }
\nc{\ug}{\mathfrak{u} } \nc{\dg}{\mathfrak{d} } \nc{\ag}{\mathfrak{a} }
\nc{\pg}{\mathfrak{p} } \nc{\sg}{\mathfrak{s} } \nc{\rca}{\mathfrak{R}}
\nc{\pca}{\mathcal{P}} \nc{\nca}{\mathcal{N}} \nc{\lca}{\mathcal{L}}
\nc{\oca}{\mathcal{O}} \nc{\mca}{\mathcal{M}} \nc{\tca}{\mathcal{T}}
\nc{\aca}{\mathcal{A}} \nc{\cca}{\mathcal{C}} \nc{\sca}{\mathcal{S}}
\nc{\bca}{\mathcal{B}}

\nc{\vp}{\varphi} \nc{\ddt}{{\small \frac{{\rm d}}{{\rm d}t}}} \nc{\im}{\mathtt{i}}

\nc{\SO}{{\mathrm SO}} \nc{\Spe}{{\mathrm Sp}} \nc{\Sl}{{\mathrm SL}}
\nc{\SU}{{\mathrm SU}} \nc{\Or}{{\mathrm O}} \nc{\U}{{\mathrm U}} \nc{\Gl}{{\mathrm
GL}} \nc{\Se}{{\mathrm S}} \nc{\Cl}{{\mathrm Cl}} \nc{\Spein}{{\mathrm Spin}}
\nc{\Pin}{{\mathrm Pin}}

\nc{\RR}{{\Bbb R}} \nc{\HH}{{\Bbb H}} \nc{\CC}{{\Bbb C}} \nc{\ZZ}{{\Bbb Z}}
\nc{\FF}{{\Bbb F}} \nc{\NN}{{\Bbb N}} \nc{\QQ}{{\Bbb Q}} \nc{\PP}{{\Bbb P}}
\nc{\G}{\mathrm{GL}_n(\RR)}
\nc{\preq}{\simeq_\RR}
\nc{\prek}{\simeq_K}

\nc{\vs}{\vspace{.2cm}} \nc{\vsp}{\vspace{1cm}} \nc{\ip}{\langle\cdot,\cdot\rangle}
\nc{\la}{\langle} \nc{\ra}{\rangle} \nc{\unm}{\frac{1}{2}} \nc{\unc}{\frac{1}{4}}
\nc{\und}{\frac{1}{16}} \nc{\no}{\vs\noindent} \nc{\lam}{\Lambda^2\ngo^*\otimes\ngo}
\nc{\tangz}{{\rm T}^{\rm Zar}} \nc{\nor}{{\sf n}}
\nc{\eigen}{(k_1<...<k_r;d_1,...,d_r)} \nc{\eigencero}{(0<k_2<...<k_r;d_1,...,d_r)}
\nc{\mum}{/\!\!/} \nc{\kir}{/\!\!/\!\!/}
\nc{\lb}{[\cdot,\cdot]}

\nc{\He}{\operatorname{Hess}} \nc{\ad}{\operatorname{ad}}
\nc{\Ad}{\operatorname{Ad}} \nc{\rank}{\operatorname{rank}}
\nc{\Irr}{\operatorname{Irr}} \nc{\End}{\operatorname{End}}
\nc{\Aut}{\operatorname{Aut}} \nc{\Inn}{\operatorname{Inn}}
\nc{\Der}{\operatorname{Der}} \nc{\Ker}{\operatorname{Ker}}
\nc{\Iso}{\operatorname{I}} \nc{\Diff}{\operatorname{Diff}}
\nc{\Lie}{\operatorname{L}} \nc{\tr}{\operatorname{tr}} \nc{\dif}{\operatorname{d}}
\nc{\sen}{\operatorname{sen}} \nc{\modu}{\operatorname{mod}}
\nc{\Ric}{\operatorname{Ric}} \nc{\Ricac}{\operatorname{Ric^{ac}}}
\nc{\Ricg}{\operatorname{Ric^{\gamma}}} \nc{\Ricc}{\operatorname{Ric^{c}}}
\nc{\Ricj}{\operatorname{Ric^{J}}}
\nc{\sym}{\operatorname{sym}} \nc{\symac}{\operatorname{sym^{ac}}}
\nc{\symc}{\operatorname{sym^{c}}} \nc{\scalar}{\operatorname{sc}}
\nc{\grad}{\operatorname{grad}} \nc{\ricci}{\operatorname{Rc}}
\nc{\ricciac}{\operatorname{ric^{ac}}} \nc{\riccic}{\operatorname{ric^{c}}}
\nc{\riccig}{\operatorname{ric^{\gamma}}} \nc{\Rin}{\operatorname{M}}
\nc{\Le}{\operatorname{L}} \nc{\tang}{\operatorname{T}}
\nc{\level}{\operatorname{level}} \nc{\rad}{\operatorname{r}}
\nc{\abel}{\operatorname{ab}}
\nc{\Pf}{\operatorname{Pf}}
\nc{\Dg}{\operatorname{Diag}}
\nc{\spa}{\operatorname{span}}

\theoremstyle{plain}
\newtheorem{theorem}{Theorem}[section]

\newtheorem{corollary}[theorem]{Corollary}

\theoremstyle{definition}
\newtheorem{definition}[theorem]{Definition}

\theoremstyle{remark}
\newtheorem{remark}[theorem]{Remark}

\newtheorem{example}[theorem]{Example}

\title[Minimal metrics]{Minimal metrics on $6$-dimensional complex nilmanifolds}

\author{Edwin Alejandro Rodr\'{\i}guez Valencia}

\address{FaMAF and CIEM, Universidad Nacional de C\'ordoba, Medina Allende s/n, 5000 C\'ordoba, Argentina}
\email{earodriguez@famaf.unc.edu.ar}

\thanks{
This research was partially supported by grants from CONICET, FONCYT and SeCyT (Universidad Nacional de C\'ordoba).}

\begin{document}
\renewcommand{\tablename}{Table}
\maketitle

\begin{abstract}
Let $(N, J)$ be a real $2n$-dimensional nilpotent Lie group endowed
with an invariant complex structure.  A left-invariant Riemannian metric
on $N$ compatible with $J$ is said to be {\it minimal}, if it minimizes the norm of the
invariant part of the Ricci tensor among all compatible metrics on $(N, J)$ with the same
scalar curvature.

In this paper, we determine all complex structures that admit a minimal compatible metric on
$6$-dimensional nilpotent Lie groups.
\end{abstract}


\section{Introduction}\label{intro}

Let $N$ be a real $2n$-dimensional nilpotent Lie group with Lie algebra $\ngo$,
whose Lie bracket will be denoted by $\mu :\ngo\times\ngo\longrightarrow\ngo$.  An {\it invariant complex structure} on $N$ is defined by a map $J:\ngo\longrightarrow\ngo$ satisfying $J^2=-I$ and the integrability condition
\begin{equation}\label{integral}
\mu(JX,JY)=\mu(X,Y)+J\mu(JX,Y)+J\mu(X,JY), \qquad \forall X,Y\in\ngo.
\end{equation}
By left translating $J$, one obtains a complex manifold $(N,J)$, as well as compact complex manifolds $(N/\Gamma,J)$ if $N$ admits cocompact discrete subgroups $\Gamma$, which are usually called \textit{nilmanifolds} and play an important role in complex geometry.

A left-invariant metric which is
{\it compatible} with $(N,J)$, also called a {\it hermitian metric}, is determined by an inner product $\ip$ on $\ngo$ such that
$$
\la JX,JY\ra=\la X,Y\ra, \quad \forall X,Y\in\ngo.
$$

A very natural evolution equation for hermitian metrics on a fixed complex manifold $(N,J)$ is given by
$$
\ddt \ip_t=-2\riccic_{\ip_t},
$$
which will be called the {\it complexified Ricci flow} (cxRF), where $\riccic_{\ip_t}:=\la\Ricc_{\ip_t}\cdot,\cdot\ra$ is the $(1,1)$-component of the Ricci tensor and $\Ricc_{\ip_t}$ the $J$-invariant part of the Ricci operator of the hermitian manifold $(N,J,\ip_t)$.

The cxRF has been studied in \cite{canonical}, where besides the uniqueness
of cxRF-\textit{solitons} up to isometry and scaling on a given $(N,J)$, the following characterizations were given:
\begin{itemize}
\item[(i)] $\ip$ is a cxRF-soliton.
\item[(ii)] $\ip$ is \textit{minimal}, that is, it minimizes the functional $\tr(\Ricc_{\ip})^2$
on the set of all compatible metrics on $(N,J)$ with the same scalar curvature.
\item[(iii)] $\Ricc_{\ip}=cI+D$ for some $c\in\RR$ and $D\in\Der(\ngo)$.
\end{itemize}

In \cite{invcomplx}, we determined which $6$-dimensional abelian (i.e. $\mu(JX,JY)=\mu(X,Y)$) complex nilmanifolds
admit a minimal metric.  In \cite{solitonSCF}, Fern\'{a}ndez-Culma gives a criterion for determining the existence of a minimal compatible metric for a geometric structure on a nilpotent Lie group, which is based on the moment map of a real reductive representation (see Section \ref{existmetric}).

Our aim in this paper is to use equivalence (iii) above and the criterion given in \cite{solitonSCF}
to classify all complex structures admitting a minimal compatible metric on $6$-dimensional nilpotent Lie groups.
In some cases, we found the minimal metrics explicitly.  A complete classification result is given in Tables
\ref{tablemin1} and \ref{tablemin2}.


\section{Preliminaries}\label{basic}

Let $\ngo$ a $2n$-dimensional real vector space, and consider the space of all skew-symmetric
algebras of dimension $2n$, which is parameterized by the vector space
$$V = \lam = \{\mu : \ngo \times \ngo \to \ngo :
\mu \ \text{bilinear and skew-symmetric}\}.$$
We now fix a map $J:\ngo \to \ngo$ such that $J^2=-I$. There is a natural linear action of Lie group
$\Gl_n(\CC):=\{g\in\Gl_{2n}(\RR): g J = J g\}$ on $V$ defined by
\begin{align}
g\cdot\mu(X, Y)=g\mu(g^{-1} X,g^{-1} Y), \quad X,Y\in\ngo, \ g\in\Gl_n(\CC), \ \mu\in V,
\end{align}
and the corresponding representation of the Lie algebra $\glg_n(\CC)$ of $\Gl_n(\CC)$ on $V$ is given
by
\begin{align}\label{rept}
\pi(\alpha)\mu = \alpha\mu(\cdot,\cdot) - \mu(\alpha\cdot,\cdot) - \mu(\cdot,\alpha\cdot),
\quad \alpha\in\glg_n(\CC), \ \mu\in V.
\end{align}
Any inner product $\ip$ on $\ngo$ determines inner products on $V$ and $\glg_n(\CC)$, both
also denoted by $\ip$, as follows:
\begin{align}\label{prodtint}
\la \mu, \lambda \ra =\sum_{i,j,k} \la \mu(e_i,e_j), e_k \ra \la \lambda(e_i,e_j), e_k \ra,
\qquad \la \alpha,\beta\ra = \tr\alpha\beta^*,
\end{align}
where $\{e_i\}$ denote an ortonormal basis of $\ngo$ and $\beta^*$ the conjugate
transpose with respect to $\ip$.

We use $\glg_n(\CC) = \ug(n)\oplus\hg(n)$ as a Cartan decomposition of $\glg_n(\CC)$, where $\ug(n)$ and
$\hg(n)$ denote the subspaces of skew-hermitian and hermitian matrices, respectively.
The set $\ag$ of all (real) diagonal $n\times n$ matrix in $\glg_n(\CC)$ is a maximal abelian
subalgebra of $\hg(n)$ and therefore determines a system of roots $\Delta\subset\ag$.
Let $\Phi$ denote the set of roots.  If $J e_{2i-1} = e_{2i}$, $i=1,\ldots,n$, then $\Phi$ is given by
\begin{align}
\Phi = \{\pm\Dg(1,-1,0,\ldots,0), \pm\Dg(1,0,-1,\ldots,0), \pm\Dg(0,1,-1,\ldots,0),\ldots\}.
\end{align}
If $\{e^1,\ldots,e^{2n}\}$ is the basis
of $\ngo^*$ dual to the basis $\{e_1, ..., e_{2n}\}$, then
\begin{align}\label{vsubijk}
\{v_{ijk}=(e^i\wedge e^j)\otimes e_k : 1\leq i < j\leq 2n, 1\leq k \leq 2n \}
\end{align}
is a basis of weight vectors of $V$ for the representation (\ref{rept}), where $v_{ijk}$ is actually
the bilinear form on $\ngo$ defined by $v_{ijk}(e_i,e_j)= - v_{ijk}(e_j,e_i)= e_k$ and zero otherwise.
The corresponding weights $\alpha_{ij}^k \in \ag$, $i<j$, are given by
\begin{align}
\pi(\alpha)v_{ijk} = (a_k - a_i - a_j)v_{ijk} = \la\alpha,\alpha_{ij}^k\ra v_{ijk}, \quad
\forall\alpha = \left[\begin{smallmatrix} a_1&&&&\\ &&&&\\
&&\ddots&&\\ &&&&\\&&&&a_n
\end{smallmatrix}\right] \in \ag.
\end{align}
If $J e_{2i-1} = e_{2i}$, it is easy to check that
$$\alpha_{ij}^k= \frac{1}{2}(E_{k,k} + E_{k\mp 1,k\mp 1} - E_{i,i} - E_{i\mp 1,i\mp 1}
- E_{j,j} - E_{j\mp 1,j\mp 1}),$$ where $E_{r,s}$ denotes the matrix whose only nonzero coefficient is
$1$ at entry $rs$.


\section{Minimal metrics on complex nilmanifolds}\label{existmetric}

Let $N$ be a real $2n$-dimensional nilpotent Lie group with Lie algebra $\ngo$, and
$J$ an invariant complex structure on $N$. A left invariant metric which is
{\it compatible} with the nilmanifold $(N,J)$, also called a {\it hermitian metric}, is determined
by an inner product $\ip$ on $\ngo$ such that
$$
\la JX,JY\ra=\la X,Y\ra, \quad \forall X,Y\in\ngo.
$$
We consider
$$
\Ricc_{\ip}:=\unm\left(\Ric_{\ip}-J\Ric_{\ip}J\right),
$$
the complexified part of the Ricci operator $\Ric_{\ip}$ of the hermitian manifold $(N,J,\ip)$, and the corresponding
$(1,1)$-component of the Ricci tensor $\riccic_{\ip}:=\la\Ricc_{\ip}\cdot,\cdot\ra$.

A compatible metric $\ip$ on $(N,J)$ is called {\it minimal} if
$$
\tr{(\Ricc_{\ip})^2}=\min \left\{ \tr{(\Ricc_{\ip'})^2} :
\scalar(\ip')=\scalar(\ip)\right\},
$$
where $\ip'$ runs over all compatible metrics on $(N,J)$ and $\scalar(\ip)=\tr\Ric_{\ip}=\tr\Ricc_{\ip}$ is the scalar curvature.  In \cite{canonical}, the following conditions on $\ip$ are proved to be equivalent to minimality:
\begin{itemize}
\item[(i)] The solution $\ip_t$ with initial value $\ip_0=\ip$ to the cxRF
$$
\ddt \ip_t=-2\riccic_{\ip_t},
$$
is self-similar, in the sense that $\ip_t=c_t\vp_t^*\ip$ for some $c_t>0$ and one-parameter group
of automorphisms $\vp_t$ of $N$.  In this case, $\ip$ is called a cxRF-{\it soliton}.

\item[(ii)] There exist a vector field $X$ on $N$ and $c\in\RR$ such that
$$\riccic_{\ip}=c\ip+L_X\ip,$$ where $L_X\ip$ denotes the usual Lie derivative.

\item[(iii)] $\Ricc_{\ip}=cI+D$ for some $c\in\RR$ and
$D\in\Der(\ngo)$.
\end{itemize}

The uniqueness up to isometry and scaling of a minimal metric on a
given $(N,J)$ was also proved in \cite{canonical}, and can be used to obtain invariants in the following way.  If $(N,J_1,\ip_1)$ and $(N,J_2,\ip_2)$ are minimal and $J_1$ is equivalent to $J_2$ (i.e. if there
exists an automorphism $\alpha$ of $\ngo$ satisfying $J_2 = \alpha J_1 \alpha^{-1}$), then they must be conjugate via an automorphism which is an isometry between $\ip_1$ and $\ip_2$.  This provides us with a lot of invariants, namely the Riemannian geometry invariants including all different kind of curvatures. In \cite{invcomplx}, we used this to
give an alternative proof of the pairwise non-isomorphism between the structures which have appeared in the
classification of abelian complex structures on $6$-dimensional nilpotent Lie algebras given in \cite{AndBrbDtt},
where condition (iii) is strongly applied.

\begin{example}\label{NOricciC}
For $t\in\RR$, consider the $3$-step nilpotent Lie algebra $\hg_{11}$ whose bracket is given by
$$
\begin{array}{lll}
\mu_{t}(e_1,e_2)= e_4, & \mu_{t}(e_1,e_3)= -e_5,\\
\mu_{t}(e_1,e_4)= (t-1)e_6, & \mu_{t}(e_2,e_3)=-t e_6.
\end{array}
$$
Let
\begin{align}\label{Jstand}
\begin{array}{lll}
J:=\left[\begin{smallmatrix} 0&-1&&&&\\ 1&0&&&&\\ &&0&-1&&\\
&&1&0&&\\ &&&&0&-1\\ &&&&1&0
\end{smallmatrix}\right], && \la e_i,e_j\ra:=\delta_{ij}.
\end{array}
\end{align}
A straightforward verification shows that $J$ is a non-abelian complex structure on $N_{\mu_{t}}$ for all
$t$ ($N_{\mu_{t}}$ is the (simply connected) nilpotent Lie group with Lie algebra $(\hg_{11}, \mu_{t})$), and $\ip$ is compatible with $(N_{\mu_{t}},J)$.  It is easy to see that $\Ricc_{\mu_{t}} = cI+D$ for some $c\in \RR$, $D\in\Der(\ngo)$ if and only if $t=0$ or $t=1$. Condition (iii) now shows that $\ip$ is not minimal for
$t>1$.
\end{example}

The problem of finding a minimal metric can be very difficult.  In \cite{solitonSCF}, Fern\'{a}ndez-Culma gives a criterion for determining the existence of a minimal compatible metric for a geometric structure on a nilpotent Lie group.  We will apply such result in the complex case.

We follow the notation of Section \ref{basic} for a fixed complex structure $J$ on $N$. Set $A=\exp(\ag)$
and consider $W$ a $A$-invariant subspace of $V$. It follows that
$W$ has a decomposition in weight spaces $$W = W_1 \oplus^{\perp}\cdots \oplus^{\perp} W_r$$
with weights $\Psi_{W} = \{\alpha_1,\ldots,\alpha_r\}$.

\begin{definition}\cite[Definition 2.18.]{solitonSCF}\label{Jnice}
We call $W$ \textit{J-nice} if $\Ricc_{\mu}\in \ag$ for all $\mu\in W$.
\end{definition}

A very useful corollary is the following
\begin{corollary}\cite[Corollary 4.7.]{distorbit}\label{critJnice}
Let $W$ be an $A$-invariant subspace of $V$. If for all $\alpha_i$ and $\alpha_j$ in $\Psi_{W}$,
$\alpha_i - \alpha_j \notin \Phi$, then $W$ is $J$-nice.
\end{corollary}

From an algebraic point of view, there is a condition on the basis of a Lie algebra that gets a subspace
$J$-nice, based on the simplicity of the corresponding set of structural constants. Namely, a basis
$\{X_1,\ldots, X_n\}$ of $\ngo$ is said to be \textit{nice} if $[X_i, X_j]$ is always a scalar multiple of
some element in the basis and two different brackets $[X_i, X_j]$, $[X_r, X_s]$ can be a nonzero
multiple of the same $X_k$ only if $\{i, j\}$ and $\{r, s\}$ are disjoint. It is easily to check that if $W$
admits a nice basis, then $W$ is $J$-nice (see \cite{riccidiag} for other application).

Let us denote by $\rca(\mu)$ the ordered set of weights related with $\mu$ to the action of
$\Gl_n(\CC)$ on $V$. It is clear that $\rca(\mu)$ is the orthogonal projection onto $\ag$ of the
weights related with $\mu$ to the action of $\Gl_{2n}(\RR)$ on $V$. We denote by
$\mathrm{U}_\mu$ the \textit{Gram matrix} of $(\rca(\mu),\ip)$, i.e.
$$\mathrm{U}_\mu(p,q) = \la \rca(\mu)_p, \rca(\mu)_q \ra$$ with
$1\leq p,q\leq \sharp\rca(\mu)$.

\begin{theorem}\cite[Theorem 2.22.]{solitonSCF}\label{thmJnice}
Let $W$ be a $J$-nice space and let $(N_\mu,J)$ be a complex nilmanifold with $\mu\in W$. \
Then $(N_\mu,J)$ admits a compatible minimal metric if and only if the equation
$$\mathrm{U}_{\mu}[x_i] = \lambda [1]$$ has a positive solution $[x_i]$  for some
$\lambda\in\RR$.
\end{theorem}

\begin{example}
By using the notation of Example \ref{NOricciC}, we will now prove that $(N_{\mu_{t}},J)$ does
admit a compatible minimal metric for all $t>1$. Let
$$W = \spa_{\RR}\{\mu_{12}^4, \mu_{13}^5, \mu_{14}^6, \mu_{23}^6\},$$
where $\mu_{ij}^k$ is defined as in (\ref{vsubijk}). Let us first see that $W$ is $J$-nice by
using Corollary \ref{critJnice}. The root set $\Phi$ of $\glg_3(\CC)$ is given by
\begin{align*}
\Phi = \{\pm\Dg(1,-1,0), \pm\Dg(1,0,-1), \pm\Dg(0,1,-1)\}.
\end{align*}
The weights of $W$ with respect to the action of $\Gl_3(\CC)$ are
{\small\begin{align*}
\{\alpha_1:=\Dg(-2,1,0), \alpha_2:=\Dg(-1,-1,1), \alpha_3:=\Dg(-1,-1,1), \alpha_4:=\Dg(-1,-1,1)\},
\end{align*}}
for all $t>1$.  Since $\alpha_i - \alpha_j \notin \Phi$, it follows that $W$ is $J$-nice.  It follows that
$$
\mathrm{U}_{\mu_{t}}=\left[\begin{smallmatrix} 5&1&1&1\\ 1&3&3&3\\ 1&3&3&3\\
1&3&3&3
\end{smallmatrix}\right].
$$
Since $X = (\frac{1}{7}, \frac{1}{7}, \frac{1}{14}, \frac{1}{14})$ is a positive solution to
the problem $\mathrm{U}_{\mu_{t}} X = [1]_4$, we conclude that $(N_{\mu_{t}},J)$ does admit a
minimal metric for all $t>1$, by Theorem \ref{thmJnice} (see Table \ref{tablemin1}).
\end{example}

\begin{example}
Consider the $2$-step nilpotent Lie algebra $(\hg_5,\mu_{st})$ given by
$$
\begin{array}{lll}
\mu_{st}(e_1,e_2)= 2 e_6, & \mu_{st}(e_1,e_3)= -e_5, & \mu_{st}(e_1,e_4)= -e_6,\\
\mu_{st}(e_2,e_3)= -e_6, & \mu_{st}(e_2,e_4)= e_5, & \mu_{st}(e_3,e_4)= 2s e_5 + 2t e_6,
\end{array}
$$
with $s\geq 0$, $t\in\RR$, $4s^2 < 1+4t$. We have that $(N_{\mu_{st}}, J)$ is a non-abelian complex
nilmanifold for all $s, t$, where $J$ is given as in (\ref{Jstand}). Let
$$W = \spa_{\RR}\{\mu_{12}^6, \mu_{13}^5, \mu_{14}^6, \mu_{23}^6, \mu_{24}^5, \mu_{34}^5,
\mu_{34}^6\}.$$
The weights of $W$ with respect to the action of $\Gl_3(\CC)$ are
{\small\begin{align*}
\{ & \alpha_1:=\Dg(-2,0,1), \alpha_2:=\Dg(-1,-1,1), \alpha_3:=\Dg(-1,-1,1), \alpha_4:=\Dg(-1,-1,1), \\
& \alpha_5:=\Dg(-1,-1,1), \alpha_6:=\Dg(0,-2,1), \alpha_7:=\Dg(0,-2,1),\},
\end{align*}}
for all $s\neq 0$, $t\neq 0$.  Since $\alpha_1 - \alpha_2 \in \Phi$, $\Phi$ as in the above example, Corollary \ref{critJnice} does not apply.  Anyway, it is straightforward to check that $\Ricc(g\cdot\mu_{st}) \in \ag$
for all $g\in A$, and so $W$ is $J$-nice.  Hence
$$
\mathrm{U}_{\mu_{st}}=\left[\begin{smallmatrix} 5&3&3&3&3&1&1\\ 3&3&3&3&3&3&3\\ 3&3&3&3&3&3&3\\
3&3&3&3&3&3&3\\ 3&3&3&3&3&3&3\\ 1&3&3&3&3&5&5\\ 1&3&3&3&3&5&5
\end{smallmatrix}\right].
$$
Since $X = (\frac{1}{12}, \frac{1}{120}, \frac{1}{40}, \frac{1}{15}, \frac{1}{15}, \frac{1}{24},
\frac{1}{24})$ is a positive solution to the problem $\mathrm{U}_{\mu_{st}} X = [1]_7$, it follows
that $(N_{\mu_{st}},J)$ does admit a minimal metric for all $(s,t)\neq (0,0)$ (analogously if $s=0$ or $t=0$),
for Theorem \ref{thmJnice}.

But if we now take $s=t=0$ then
$$W = \spa_{\RR}\{\mu_{12}^6, \mu_{13}^5, \mu_{14}^6, \mu_{23}^6, \mu_{24}^5\}, \qquad
\mathrm{U}_{\mu}=\left[\begin{smallmatrix} 5&3&3&3&3\\ 3&3&3&3&3\\ 3&3&3&3&3\\
3&3&3&3&3\\ 3&3&3&3&3
\end{smallmatrix}\right].
$$
Any solution to the equation $\mathrm{U}_{\mu} X = \lambda[1]_5$ is of the form $(0, \frac{1}{3}-a-b-c, a, b, c)$,
and therefore $(N_\mu, J)$ does not admit a minimal metric by Theorem \ref{thmJnice}. In summary,
$(N_{\mu_{st}},J)$ does admit a minimal metric if and only if $s\neq 0$ or $t\neq 0$ (see Table \ref{tablemin1}).
\end{example}

\begin{example}
We consider the $4$-step nilpotent Lie algebra $(\hg_{26}^{+},\mu)$ defined by
$$
\begin{array}{lll}
\mu(e_1,e_2)= e_5, & \mu(e_1,e_3)= \pm e_6, & \mu(e_1,e_5)= -e_3,\\
\mu(e_2,e_4)= \pm e_6, & \mu(e_2,e_5)= -e_4.
\end{array}
$$
Therefore, $(N_\mu, J)$ is a complex nilmanifold (see (\ref{Jstand})). Let
$$W = \spa_{\RR}\{\mu_{12}^5, \mu_{13}^6, \mu_{15}^3, \mu_{24}^6, \mu_{25}^4\}.$$
Note that $W$ is nice, and, in consequence, it is $J$-nice.
The weights of $W$ with respect to the action of $\Gl_3(\CC)$ are
{\small\begin{align*}
\{ & \alpha_1:=\Dg(-2,0,1), \alpha_2:=\Dg(-1,-1,1), \alpha_3:=\Dg(-1,1,-1), \alpha_4:=\Dg(-1,-1,1), \\
& \alpha_5:=\Dg(-1,1,-1)\}.
\end{align*}}
Thus
$$
\mathrm{U}_{\mu}=\left[\begin{smallmatrix} 5&3&1&3&1\\ 3&3&-1&3&-1\\ 1&-1&3&-1&3\\
3&3&-1&3&-1\\ 1&-1&3&-1&3
\end{smallmatrix}\right].
$$
Any solution to $\mathrm{U}_{\mu} X = \lambda[1]_5$ is of the form $(-2, 3-a, 2-b, a, b)$,
thus $(N_\mu, J)$ does not admit a minimal metric by Theorem \ref{thmJnice} (see Table \ref{tablemin2}).
\end{example}


\section{Classification of minimal metrics on $6$-dimensional nilpotent Lie groups}\label{table}

In this section, we use the classification of all complex structures on $6$-dimensional nilpotent Lie groups
given in \cite{nilstruc}, to determine those admitting a minimal hermitian metric.  Recall that in
\cite{invcomplx}, we was analyzed the abelian case and therefore we do not study it here.

Next we illustrate how to rewrite the complex structure equations appearing in \cite{nilstruc}
on the Lie algebra $\hg_2$, in a way that the complex structure $J$ be fixed and varies the bracket.
Let $J e^1 = e^2$, $J e^3 = e^4$ and $J e^5 = e^6$ ($J$ view in the dual $\hg_2^{\ast}$, recall that
$(J\alpha)(x)=-\alpha(Jx)$ for all $\alpha\in\Lambda^2\hg_2^*$). With respect
to the basis
$$\{\omega^1:= e^1-i J e^1, \ \omega^2:= e^3-i J e^3, \
\omega^3:= e^5-i J e^5\},$$ the complex structure equations are \ $d\omega^1 = d\omega^2 = 0,
d\omega^3 = \omega^{12} + \omega^{1\overline{1}} + \omega^{1\overline{2}} + D
\omega^{2\overline{2}}$, with $D \in \CC$ and $\mathfrak{Im} D >0$. Here $\omega^{jk}$ (resp.
$\omega^{j\overline{k}}$) means the wedge product $\omega^j\wedge \omega^k$ (resp.
$\omega^j\wedge \omega^{\overline{k}}$), where $\omega^{\overline{k}}$ indicates the complex
conjugation of $\omega^k$. Let $D=\mathfrak{Re} D + i \mathfrak{Im} D = t + i s$, $s>0$.  It
follows that
\begin{align*}
& d e^1 = d e^2 = d e^3 = d e^4 = 0,\\
& d e^5 - i d e^6 = 2i e^1\wedge e^2 + 2 e^1\wedge e^3 - 2i e^2\wedge e^3 +
2(i t - s) e^3\wedge e^4.
\end{align*}
Therefore,
\begin{align*}
& d e^1 = d e^2 = d e^3 = d e^4 = 0,\\
& d e^5 = 2 e^1\wedge e^3 - 2s e^3\wedge e^4,\\
& d e^6 = -2 e^1\wedge e^2 + 2 e^2\wedge e^3 -2 t e^3\wedge e^4.
\end{align*}
Recall that \ $d e^k = \sum\limits_{i,j}(-c_{ij}^k) e^i\wedge e^j \Leftrightarrow [e_i, e_j] =
\sum\limits_k c_{ij}^k e_k$. Hence,
\begin{align*}
& [e_1, e_2] = e_6, \ [e_1, e_3] =-e_5,\\
& [e_2, e_3] = -e_6, \ [e_3, e_4] = s e_5 + t e_6.
\end{align*}

By arguing as above for each item in \cite[Table 1, 2]{nilstruc}, and applying Theorem \ref{thmJnice},
with $J$ as in (\ref{Jstand}), we can now formulate our main result.  Let $N_4^t$, $N_5^{s,t}$ and
$N_{26}^{+}$ denote the nilpotent Lie groups with Lie algebras $(\hg_4, \lb_t)$, $(\hg_5,\lb_{s,t})$ and
$(\hg_{26}^{+},\lb_{\pm})$, respectively (see Tables \ref{tablemin1} and \ref{tablemin2}).

\begin{theorem}
Let $(N,J)$ be a $6$-dimensional complex nilmanifold.  Then $(N,J)$ admits a minimal metric if and only
if $(N,J)$ is not holomorphically isomorphic to one of the following: \ $(N_4, J)$ \ (abelian),
$\left(N_4^{1/4},J\right)$, $\left(N_5^{0,0},J\right)$ and $\left(N_{26}^{+},J^{\pm}\right)$.
\end{theorem}

\begin{remark}
In the notation used in \cite{nilstruc}, the five exceptions above correspond to the following complex structures:
\begin{align*}
\hg_4: \quad & d\omega^1 = d\omega^2 = 0, \
d\omega^3 = \omega^{1\overline{1}} + \omega^{1\overline{2}} + (1/4)
\omega^{2\overline{2}} \quad (\mbox{abelian}).\\
& d\omega^1 = d\omega^2 = 0, \
d\omega^3 = \omega^{12} + \omega^{1\overline{1}} + \omega^{1\overline{2}} + (1/4)
\omega^{2\overline{2}}.\\
\hg_5: \quad & d\omega^1 = d\omega^2 = 0, \
d\omega^3 = \omega^{12} + \omega^{1\overline{1}}.\\
\hg_{26}^{+}: \quad & d\omega^1 = 0, \ d\omega^2 = \omega^{13} + \omega^{1\overline{3}}, \
d\omega^3 = i\omega^{1\overline{1}} \pm i(\omega^{1\overline{2}}-\omega^{2\overline{1}}).
\end{align*}
\end{remark}

In Tables \ref{tablemin1} and \ref{tablemin2}, we are given explicitly the Lie algebras $\ngo$, indicating
the condition under which $(N,J)$ admits a minimal compatible metric in the third column.  In the last column, we
added the condition under which the canonical metric $\ip$ (see (\ref{Jstand})) is minimal, that is, the cases in which $\Ricc_{\ip} = cI+D$ for some $c\in \RR$, $D\in\Der(\ngo)$.
\footnotesize{
\begin{table}
\centering
\renewcommand{\arraystretch}{1.6}
\begin{tabular}{|c|l|c|c|}\hline
$\mathbf{\ngo}$ & \textbf{Bracket} & \textbf{Existence} & $\ip$ \textbf{minimal}\\ \hline\hline
$\hg_2$ & $[e_1,e_2]=e_6$, \ $[e_1,e_3]=-e_5$, \ $[e_2,e_3]=-e_6$, & Yes & No\\
& $[e_3,e_4]=s e_5 + t e_6$; \  $s>0$, $t\in\RR$. & &\\ \hline
$\hg_4$ & $[e_1,e_2]=e_6$, \ $[e_1,e_3]=-e_5$, \ $[e_2,e_3]=-e_6$, &  $t\neq\frac{1}{4}$ & $t=-1$\\
& $[e_3,e_4]= t e_6$; \  $t\in\RR-\{0\}$. & &\\ \hline
$\hg_5$ & $[e_1,e_3]=-e_5$, \ $[e_1,e_4]=-e_6$, &  Yes & Yes\\
& $[e_2,e_3]=-e_6$, \ $[e_2,e_4]=e_5$. &  &\\ \cline{2-4}
& $[e_1,e_2]=2 e_6$, \ $[e_1,e_3]=-e_5$, \ $[e_1,e_4]=-e_6$, & $s\neq 0$ or $t\neq 0$ & $s^2+t^2=1$\\
& $[e_2,e_3]=-e_6$, \ $[e_2,e_4]=e_5$, \ $[e_3,e_4]=2s e_5 + 2t e_6$, & &\\
& $s\geq 0$, \ $t\in\RR$, \ $4s^2 < 1+4t$. & &\\ \cline{2-4}
& $[e_1,e_2]=2 e_6$, \ $[e_1,e_3]=-(t+1)e_5$, \ $[e_1,e_4]=(t-1)e_6$, & Yes & No\\
& $[e_2,e_3]=-(t+1) e_6$, \ $[e_2,e_4]=(1-t)e_5$, \ $[e_3,e_4]=2s e_5$, & &\\
& with $(s,t)$ satisfying one of: & &\\
& $\bullet$ \ $0<t^2<\frac{1}{2}$, \ \ $0\leq s < \frac{t^2}{2}$ & &\\
& $\bullet$ \ $\frac{1}{2}\leq t^2<1$, \ $0\leq s < \frac{1-t^2}{2}$ & &\\
& $\bullet$ \ $t^2>1$, \ $0\leq s < \frac{t^2-1}{2}$ & &\\ \hline
$\hg_6$ & $[e_1,e_2]=e_6$, \ $[e_1,e_3]=-e_5$, \ $[e_2,e_3]=-e_6$. & Yes & No\\ \hline
$\hg_7$ & $[e_1,e_2]=e_4$, \ $[e_1,e_3]=-e_5$, \ $[e_2,e_3]=-e_6$. & Yes & Yes\\ \hline
$\hg_{10}$ & $[e_1,e_2]=e_4$, \ $[e_2,e_3]=-e_6$, \ $[e_2,e_4]=e_5$. & Yes & Yes\\ \hline
$\hg_{11}$ & $[e_1,e_2]=e_4$, \ $[e_1,e_3]=-t e_5$, \ $[e_2,e_3]=-e_6$, &  Yes & No\\
& $[e_2,e_4]=(1-t) e_5$; \ $t<1$, $t\neq 0$. &  &\\ \cline{2-4}
& $[e_1,e_2]=e_4$, \ $[e_1,e_3]=-e_5$, \ $[e_1,e_4]=(t-1)e_6$, & Yes & No\\
& $[e_2,e_3]=-t e_6$; \ $t>1$. &  &\\ \hline
$\hg_{12}$ & $[e_1,e_2]=2e_4$, \ $[e_1,e_3]=-(s+1-\alpha)e_5 + t e_6$, & Yes &
$\left(s-\frac{1}{2}\right)^2+t^2=\frac{1}{4}$\\
& $[e_1,e_4]=t e_5 + (s-1+\alpha)e_6$, \ $[e_2,e_3]=-t e_5 - (s+1+\alpha) e_6$, & &\\
& $[e_2,e_4]=-(s-1-\alpha)e_5 + t e_6$; \ $s, t \in\RR$, \ $t\neq 0$, & &\\
& with $\alpha:=\sqrt{(s-1)^2+t^2}$. & &\\ \hline
$\hg_{13}$ & $[e_1,e_2]=2e_4$, \ $[e_1,e_3]=-(s+1-c)e_5 + t e_6$, &  Yes & $c^2+s^2+t^2=1$\\
& $[e_1,e_4]=t e_5 + (s-1+c)e_6$, \ $[e_2,e_3]=-t e_5 - (s+1+c) e_6$, & &\\
& $[e_2,e_4]=-(s-1-c)e_5 + t e_6$; \ $s, t \in\RR$, \ $c\in\RR^{\geq 0}$, & &\\
& with $c,\beta$ satisfying: let $\alpha:=\sqrt{(s-1)^2+t^2}$, $\beta:=\sqrt{s^2+t^2}$, & &\\
& $c\neq \alpha$, \ $(c,\beta)\neq(0,1)$, \ $c^4-2(\beta^2+1)c^2+(\beta^2-1)^2 < 0$. & &\\ \hline
$\hg_{14}$ & $[e_1,e_2]=2e_4$, \ $[e_1,e_3]=-(s+1-c)e_5 + t e_6$, &  Yes & $c^2+s^2+t^2=1$\\
& $[e_1,e_4]=t e_5 + (s-1+c)e_6$, \ $[e_2,e_3]=-t e_5 - (s+1+c) e_6$, & &\\
& $[e_2,e_4]=-(s-1-c)e_5 + t e_6$; \ $s, t \in\RR$, \ $c\in\RR^{\geq 0}$, & &\\
& with $c,\beta$ satisfying: let $\alpha:=\sqrt{(s-1)^2+t^2}$, $\beta:=\sqrt{s^2+t^2}$, & &\\
& $c\neq \alpha$, \ $(c,\beta)\neq(0,1)$, \ $c^4-2(\beta^2+1)c^2+(\beta^2-1)^2 = 0$. & &\\ \hline
\end{tabular}
\end{table}}

\begin{table}
\centering
\renewcommand{\arraystretch}{1.6}
\begin{tabular}{|c|l|c|c|}\hline
$\mathbf{\ngo}$ & \textbf{Bracket} & \textbf{Existence} & $\ip$ \textbf{minimal}\\ \hline\hline
$\hg_{15}$ & $[e_1,e_2]=2e_4$, \ $[e_1,e_3]=-(s+1-c)e_5 + t e_6$, &  Yes & $c^2+s^2+t^2=1$\\
& $[e_1,e_4]=t e_5 + (s-1+c)e_6$, \ $[e_2,e_3]=-t e_5 - (s+1+c) e_6$, & &\\
& $[e_2,e_4]=-(s-1-c)e_5 + t e_6$; \ $s, t \in\RR$, \ $c\in\RR^{\geq 0}$, & &\\
& with $c,\beta$ satisfying: let $\alpha:=\sqrt{(s-1)^2+t^2}$, $\beta:=\sqrt{s^2+t^2}$, & &\\
& $c\neq \alpha$, \ $(c,\beta)\neq(0,1)$, \ $c^4-2(\beta^2+1)c^2+(\beta^2-1)^2 > 0$. & &\\ \hline
$\hg_{16}$ & $[e_1,e_2]=2e_4$, \ $[e_1,e_3]=-(s+1)e_5 + t e_6$, &  Yes & $s^2+t^2=1$\\
& $[e_1,e_4]=t e_5 + (s-1)e_6$, \ $[e_2,e_3]=-t e_5 - (s+1) e_6$, & &\\
& $[e_2,e_4]=(1-s)e_5 + t e_6$; \ $s, t \in\RR$, \ $s^2+t^2=1$, $(s,t)\neq(1,0)$. & &\\ \hline
\end{tabular}
\vspace{0.3cm}
\caption{Non-abelian Nilpotent complex structures.}\label{tablemin1}\vspace{-0.3cm}
\end{table}

\begin{table}
\centering
\renewcommand{\arraystretch}{1.6}
\begin{tabular}{|c|l|c|c|}\hline
$\mathbf{\ngo}$ & \textbf{Bracket} & \textbf{Existence} & $\ip$ \textbf{minimal}\\ \hline\hline
$\hg_{19}^{-}$ & $[e_1,e_3]=\pm e_6$, \ $[e_1,e_5]=-e_3$, & Yes & Yes\\
& $[e_2,e_4]=\pm e_6$, \ $[e_2,e_5]=-e_4$. & &\\ \hline
$\hg_{26}^{+}$ & $[e_1,e_2]= e_5$, \ $[e_1,e_3]=\pm e_6$, \ $[e_1,e_5]= -e_3$, & No & ------\\
& $[e_2,e_4]=\pm e_6$, \ $[e_2,e_5]=-e_4$. & &\\ \hline
\end{tabular}
\vspace{0.3cm}
\caption{Non-nilpotent complex structures.}\label{tablemin2}
\end{table}


\end{document}